\documentclass[11pt,notitlepage]{article}
\usepackage{latexsym,amsfonts,amssymb,amsmath,amsthm}
\usepackage{verbatim}
\newcommand{\lskew}{\!{<}\hspace{-.55em}(\hspace{.2em}}
\newcommand{\rskew}{{>}\hspace{-.65em})\hspace{.2em}}
\newcommand{\newatop}[2]{\genfrac{}{}{0pt}{}{#1}{#2}}
\newcommand{\newstack}[2]{\raisebox{-.3ex}{$\stackrel{#1}{\scriptstyle \,#2\,}$}}
\newcommand{\Q}{\mathbb Q}
\newcommand{\C}{\mathbb C}
\newcommand{\Sym}{\mathfrak S}

\newcommand\Ell[2]{\ell(#1 | #2)}
\newcommand\Det[1]{\mathrm{Det}\,#1}

\swapnumbers
\theoremstyle{plain}		
	\newtheorem{mytheo}{Theorem}[section]
	
	\newtheorem{myprop}[mytheo]{Proposition}
	\newtheorem{mycoro}[mytheo]{Corollary}

\theoremstyle{definition}	
	\newtheorem{mydefi}{Definition}[section]

\theoremstyle{remark}	
	\newtheorem*{mynota}{Notation}	
	\newtheorem*{mylemm}{Lemma}
{\swapnumbers	}
	\newtheorem*{myrema}{Remark}
	\newtheorem*{myremas}{Remarks}

\title{Flag Varieties for the Yangian $Y(\mathfrak{gl}_n)$}
\author{Aaron Lauve}
\date{January 3, 2006}

\begin{document}
\maketitle

\begin{abstract}
It is well-known that the Yangian $Y_n$ over $\mathfrak{gl}_n$ somewhat resembles the universal enveloping algebra for $\mathfrak{gl}_n$. In this work, we show it also possesses some features of the ring of regular functions on $\mathrm{GL}_n$. In particular, we use the theory of quasideterminants to construct noncommutative flags associated to the ring $Y_n[[u^{-1}]]$. In so doing, a class $\mathcal{F}\ell(\gamma)$ of comodule algebras for $Y_n$ (viewed as $\C[\mathrm{GL}_n]$) is revealed which, as in the classical case, contain the irreducible highest-weight modules for $Y_n$ (viewed as $U(\mathfrak{gl}_n)$). In the course of defining the rings $\mathcal{F}\ell(\gamma)$, connections to the new parabolic presentations of $Y_n$ given by Brundan and Kleshchev (2005) are uncovered.
\end{abstract}

%
\section*{Introduction}
The Yangians were introduced twenty years ago in the study of the Yang-Baxter equation 
(independently by Drinfeld \cite{Dri:1} and Jimbo \cite{Jim:1}), and in relation to the inverse scattering method (in the St.-Petersburg school, Faddeev, Takhtajan, et al \cite{TakFad:1, KulSkl:1}). An excellent and detailed account of the history and applications of the Yangians appears in Molev's survey article \cite{Mol:1}. 

Each Yangian $Y(\mathfrak{a})$ (there is one for each simple finite-dimensional Lie algebra $\mathfrak{a}$, and also for $\mathfrak{a}=\mathfrak{gl}_n$) is a deformation of the universal enveloping algebra $U(\mathfrak{a})$ for the polynomial current Lie algebra $\mathfrak{a}$ over $\C$. The deformation is such that $U(\mathfrak{a}[x])$ exists as a subalgebra of $Y(\mathfrak{a})$. For the remainder of the paper, we focus on $\mathfrak{a} = \mathfrak{gl}_n$ and write $Y_n$ for the associated Yangian. 

After the preceding paragraph, it is not surprising to learn that the representations of $\mathfrak{gl}_n$ play an important role in the representation theory of $Y_n$. While this theory will make an appearance in the sequel, it is not the focus of our efforts. Our main goal is to introduce some elementary \emph{co-representations} in a novel way. 

\subsection*{Summary of results} 
In this paper, we use the quasideterminant of Gelfand and Retakh \cite{GelRet:1,GelRet:3} to 
produce a class of $Y_n$-comodule algebras which may be viewed as (coordinate rings of) flag varieties for the Yangian. We show that, analogous to the classic setting, these algebras comprise irreducible highest-weight modules for $Y_n$. Section \ref{sec:comm-flags} reviews the classic construction of flag varieties and their homogeneous coordinate rings. After introducing the Yangian and its determinant in Section \ref{sec:yangian}, we use the theory of quasideterminants to discover the main object of study in Section \ref{sec:NC-flags}. We conclude in Section \ref{sec:main-results} with the main results stated above.

\subsection*{Notation}
We fix some notations and conventions used in the sequel: 
\begin{itemize}
\item[] 
Given a positive integer $n$, say $\gamma\models n$ or $\gamma$ is a \emph{composition} of $n$ if $\gamma$ of $n$ is a sequence of \emph{positive} integers summing to $n$.
\item[] 
By $[n]$ we mean the set $\{1,2,\ldots n\}$; by $[n]^k$ we mean the set of all $k$-tuples chosen from $[n]$; and by ${\scriptstyle \binom{[n]}{k}}$ we mean the set of all subsets of $[n]$ of size $k$. 
\item[] 
For two integers $m,n$ and two subsets $I\subseteq[m]$ and $J\subseteq[n]$ we define two matrices associated to an $m\times n$ matrix $A$. By $A^{I,J}$ we mean the
matrix obtained by deleting rows $I$ and columns $J$ from $A$. By $A_{I,J}$ we mean the matrix obtained by keeping only rows $I$ and columns $J$ of $A$. 
\item[] 
With slight abuse of the just-defined notation, $A^{ij}$ will represent the $(m-1)\times(n-1)$ minor of $A$ obtained by deleting row $i$ and column $j$. \emph{Also, $A_I$
will denote the square matrix obtained from $A$ by taking column-set $I$ and
row-set the first $|I|$ rows of $A$.}
\item[] 
For a misordered set of distinct integers $I=(i_1,i_2,\ldots,i_m)$ we denote by $\ell(I)$ the \emph{length} of the permutation represented by $I$, i.e. the minimal number of adjacent
swaps necessary to order $I$. We write $\ell(\sigma)$ for $\ell(\sigma1,\sigma2,\cdots,\sigma m)$.
\item[] 
All of our division rings contain $\Q$, all rings and algebras are unital.
\end{itemize}

\def\OtherRings{For a treatment of flags over any commutative ring of characteristic $p$ not dividing $n!$, see \cite{Tow:1}.}
\def\PluckerMap{For a geometric proof of the ``embedding'' part, see \cite{GriHar:1}; for an algebraic proof, see \cite{Ful:1}.}
\section{Review of Classical Setting}\label{sec:comm-flags}
We recall some classical properties of flags over $\C$ which will be mimiced for the Yangian in subsequent sections.\footnote{\OtherRings} 

\subsection{Flags}\label{sec:cmplx-flags}
Fix a vector space $V\simeq\C^n$ and a composition $\gamma=(\gamma_1,\ldots,\gamma_{r})$ of $n$.

\begin{mydefi} A \emph{flag} $\Phi$ of shape $\gamma$ is an increasing chain of subspaces of $V$, 
$$
\Phi: (0)=W_0\subsetneq W_1\subsetneq\cdots\subsetneq W_r = V\,,
$$ 
satisfying $\mathrm{dim}_{\C}\left( W_i/W_{i-1}\right) = \gamma_i$. For fixed $V$ and $\gamma$, we let $F\ell(\gamma)$ denote the collection of all flags in $V$ of shape $\gamma$. 
\end{mydefi}

\begin{mynota} Two important special cases are when $\gamma=(1^n)$ and $\gamma=(d,n-d)$. The former is the collection of \emph{full flags}, $\mathrm{dim}\,W_i=i,\,1\leq i\leq n$; the latter is the \emph{Grassmannian,} i.e. the collection of $d$-dimensional subspaces of $V$. Write $F\ell(n)$ and $Gr(d,n)$ in the respective cases.
\end{mynota}

If we fix a basis $\mathbf B^*=(f_1,\ldots,f_n)$ for $V^*$, we may represent a flag $\Phi$ as a matrix as follows. (i) Choose a basis $(w_1,\ldots,w_{\gamma_1})$ for $W_1$. (ii) Extend this to a basis $(w_1,\ldots,w_{\gamma_1},w_{\gamma_1+1},\ldots,w_{\gamma_1+\gamma_2})$ for $W_2$. (iii) Repeat until you have completed the sequence to a basis $\mathbf w=(w_1,\ldots,w_{|\gamma|})$ of $V$. (iv) Define the matrix $A=A(\mathbf w)=(a_{ij})$ by putting $a_{ij} = f_j(w_i)$. Then $A$ is the collection of row vectors $[w_1 |w_2 |\cdots |w_{n}]^T$, with the $w_i$ coordinatized by $\mathbf B$.

\begin{mylemm} Fix $\mathbf{B}$, $\Phi, \mathbf w$, and $A(\mathbf{w})$ as above. A set $\mathbf{w}'$ is another basis for $\Phi$ if and only if $A(\mathbf{w}') = g\cdot A(\mathbf{w})$ for some $g\in\mathrm{GL}_n(\C)$ of the form appearing in Figure \ref{fig:block-lower-triang}.
\begin{figure}[!hbt]
$$
\begin{array}{r@{}l}
{\begin{array}{c}\gamma_1 \vphantom{\rule[-8pt]{.5pt}{20pt}} \\ 
\gamma_2 \vphantom{\rule[-16pt]{.5pt}{35pt}} \\ 
\vdots \\ 
\gamma_r \vphantom{\rule[-9pt]{.5pt}{20pt}} 
\end{array}} &
{\left[\begin{array}{c@{}c@{}c@{}c}
\framebox[20pt][c]{$\vphantom{\rule[-5pt]{.5pt}{15pt}}g_1$} & 0 & 0 & 0 \\
 \ast & \framebox[30pt][c]{$\vphantom{\rule[-10pt]{.5pt}{25pt}}g_2$} & 0 & 0 \\
 \ast & \ast & \,\,\ddots\,\, & 0 \\
 \ast & \ast & \ast & \framebox[20pt][c]{$\vphantom{\rule[-5pt]{.5pt}{15pt}}g_r$}
\end{array}\right]}
\end{array}
$$
\caption{A lower block-triangular matrix, with $g_i\in\mathrm{GL}_{\gamma_i}(\C)$.}\label{fig:block-lower-triang}
\end{figure}
\end{mylemm}

For fixed $\gamma$, the collection of such $g\in\mathrm{GL}_n$ is a \emph{parabolic subgroup} we shall denote by $\mathrm{P}_{\gamma}$. Toward the goal of coordinatizing our flags, we replace the above definition with a new one.

\begin{mydefi} Given a composition $\gamma\models n$, we identify $Fl(\gamma)$ with the right cosets $\mathrm{P}_{\gamma}\, \backslash\, \mathrm{GL}_n(\C)$.
\end{mydefi}

\subsection{Determinants \& Coordinates}
Given a composition $\gamma \models n$, let $d_i$ denote the sum $\gamma_1 + \cdots + \gamma_i$.
Consider the map $\eta_i:Fl(\gamma) \rightarrow \mathbb P(\C^{\binom{n}{d_i}}) $ which sends $A(\Phi, \mathbf{w})$ to the $\binom{n}{d_i}$-tuple of all minors one can possibly make from the first $d_i$ rows of $A$ (not repeating columns, and taking chosen columns in order). This tuple is rightly viewed as projective coordinates because (i) it misses $0$, and (ii) it's only defined up to nonzero scalars: 
\begin{enumerate}
\item[(i)] As $A_{[d_i],[n]}$ has full rank for all $i$, there must exist one minor of size $d_i$ which is nonzero.
\item[(ii)] We need $\eta_i(gA) \equiv \eta_i(A)$ for $g\in \mathrm{P}_{\gamma}$, but the former (cf. the depiction of $g$ in Figure \ref{fig:block-lower-triang}) equals $(\prod_{j\leq i}\det g_j) \cdot \eta_i(A)$.
\end{enumerate}

We put all of these maps together to build a map $\eta: Fl(\gamma)\rightarrow \mathbb P(\gamma) := \mathbb P^{\binom{n}{d_1}-1}\times \cdots \times \mathbb P^{\binom{n}{d_{r-1}}-1}$. This map is called the \emph{Pl\"ucker embedding.}\footnote{\PluckerMap} Note that we stop at $i=r-1$. There is nothing to gain by including the final factor ($\mathbb P^0$). 

Represent a point $\pi\in\mathbb P(\gamma)$ by its coordinates $\pi=(p_I)_{I\in \binom{[n]}{\Vert \gamma \Vert}}$. When $\pi$ belongs to the image of $\eta$ ---i.e. when $\exists A\in\mathrm{GL}_n(\C)$ with (writing $|I|=d$) $p_I = \det A_{[d],I}$ for all $I\in\binom{[n]}{\Vert \gamma \Vert}$---we say the $\{p_I\}$ are the \emph{Pl\"ucker coordinates} of $A$.

The image of $\eta$ is particularly nice, it is given by quadratic relations among the coordinates $p_I$.

\begin{mytheo}\label{thm:comm-coords-characterize} For a given $\gamma\models n$ and $\pi\in\mathbb P(\gamma)$, $\pi$ belongs to the image of $\eta$ if and only if for all subsets $I=\{i_1,\ldots,i_{d-s}\}$ and $J=\{j_1,\ldots,j_{e+s}\}$ of $[n]$, for all $1\leq s$ and $d,e\in \Vert \gamma \Vert$ satisfying $d\leq e$, $\pi$ satisfies the \emph{Young symmetry relations} $(\mathcal Y_{I,J})_{(s)}$:
\begin{equation}\label{eq:cplx-YS}
0 = \sum_{{\Lambda\subseteq J},\,{|\Lambda|=s}}(-1)^{\Ell{\Lambda}{J\setminus\Lambda}} p_{I|\Lambda} p_{J\setminus\Lambda}\,.
\end{equation}
\end{mytheo}

\begin{myrema}
Here, we have extended the definition of $p_K$ from $K\in\binom{[n]}{d}$ to $K\in[n]^d$ at the expense of adding the obvious \emph{alternating relations} $(\mathcal A_K)$: 
\begin{equation}\label{eq:cplx-alternating}
p_{\sigma K} = (-1)^{\ell(\sigma)} p_{K} \qquad (\forall K\in[n]^d,\, \forall\sigma\in\Sym_d)
\end{equation}
\end{myrema}

Informed of the previous theorem, we make the following 
\begin{mydefi}\label{def:cplx-flag-algebra} The \emph{flag algebra} $\mathcal{F}\ell(\gamma)$, i.e., the multihomogeneous coordinate ring of the flag variety $Fl(\gamma)$, is the commutative $\C$-algebra with generators $\left\{f_I\mid I\in [n]^{\Vert\gamma\Vert}\right\}$ and relations $(\mathcal Y_{I,J})$ and $(\mathcal A_K)$ for allowable choices $I,J,K$.
\end{mydefi}

\subsection{Comodules}\label{sec:comm-comodules}
For fixed $\gamma\models n$, $\mathrm{GL}_n$ acts transitively on $Fl(\gamma)$ by right multiplication: $A\mapsto A' =A\cdot g$, a representative of a (possibly) different coset in 
$\mathrm{P}_{\gamma}\,\backslash\,\mathrm{GL}_n$. This representation will yield a co-representation presently. 
View $\mathrm{GL}_n$ as a variety, i.e. the open set in $\C^{n^2}$ described by the nonvanishing of the function $\det X$, where $X=(x_{ij})$ is the matrix of coordinate functions for the affine space $\C^{n^2}$. Its ring of regular functions $\C[\mathrm{GL}_n]$ is the commutative algebra generated by the $n^2+1$ generators $X=(x_{ij})$ and $y$ and the relation $\det X \cdot y - 1 = 0$. Recall that $\C[\mathrm{GL}_n]$ is a Hopf algebra with co-structure given by $\Delta(x_{ik}) = \sum_{j} x_{ij}\otimes x_{jk}$ and $\varepsilon(x_{ik}) = \delta_{ik}$. 

The mapping $Fl(\gamma)\times \mathrm{GL}_n \rightarrow Fl(\gamma)$ described above becomes an algebraic map between two varieties; we deduce the existence of an equal but opposite mapping between their rings of coordinate functions. 

\begin{myprop}\label{comm:comodule-algebra} Let $\mathcal{F}\ell_1$ be the span of the generators $\{f_I\}$ of the flag algebra. The vector space map $\rho: \mathcal{F}\ell_1 \rightarrow \mathcal{F}\ell_1 \otimes \C[\mathrm{GL}_n]$ given by
\begin{equation}\label{eq:comodule-map}
\rho(f_I) = \sum_{{J\subseteq [n]},\,{|J|=|I|}} f_J \otimes \det X_{J,I}
\end{equation}
may be extended multiplicatively to a well-defined algebra map from $\mathcal{F}\ell(\gamma)$ to $\mathcal{F}\ell(\gamma) \otimes \C[\mathrm{GL}_n]$. Moreover, $\rho$ gives $\mathcal{F}\ell(\gamma)$ the structure of right $\C[\mathrm{GL}_n]$-comodule algebra, i.e. $(\rho\otimes 1)\circ \rho = (1\otimes \Delta)\circ \rho$.
\end{myprop}

The next important result, essentially coming for free after the Pl\"ucker embedding, is the

\begin{myprop}\label{thm:comm-flag-is-subalg} The algebra $\mathcal{F}\ell(\gamma)$ is isomorphic to the subalgebra of $\C[\mathrm{GL}_n]$ generated by the minors $\big\{\det X_{[d],I} \mid I\in\binom{[n]}{\Vert \gamma \Vert} \big\}$. 
\end{myprop}

\subsection{Modules}\label{sec:comm-modules}
The universal enveloping algebra $U(\mathfrak{gl}_n)$ has generators $E_{ij}$ $(1\leq i,j\leq n)$ and relations $[E_{ij},E_{kl}] = \delta_{jk}E_{il} - \delta_{li}E_{kj}$. The generators $E_{ii}$ play a special role. 
A vector $v$ in a module $M$ for $U(\mathfrak{gl}_n)$ is called a \emph{weight vector} if there are scalars $\lambda_i$ ($1\leq i\leq n$) such that $E_{ii}\cdot v = \lambda_{i} v$; it is called a \emph{highest-weight} vector if furthermore $E_{ij} \cdot v =0$ if $i<j$. Call $M$ a \emph{highest-weight module} if $U(\mathfrak{gl}_n)\cdot v = M$. In case $\lambda=(\lambda_1, \cdots, \lambda_n)$ is a partition, i.e.  $\lambda_i\geq \lambda_{i+1}\geq 0\, (\forall i)$, write $M^{\lambda}$ to denote this special module. 

The finite dimensional, polynomial, irreducible modules for $U(\mathfrak{gl}_n)$ are understood; they are precisely the highest-weight modules $M^{\lambda}$ ($\lambda$ running over all partitions with at most $n$ parts). We next recall how $\mathcal{F}\ell(\gamma)$ comprises a sum of such $M^{\lambda}$. 

Define an action of $U(\mathfrak{gl}_n)$ on $\mathcal{F}\ell(\gamma)_1$, by
$$
E_{ab} \cdot f_I = \left\{\begin{array}{ll}
0 & \hbox{if }b\not\in I \\
f_{i_1\cdots a\cdots i_d} & \hbox{otherwise, replacing }b\hbox{ with }a.
\end{array}\right.
$$
Extend the action to all of $\mathcal{F}\ell(\gamma)$ by letting $E_{ab}$ act as a derivation. This action is well-defined, i.e. provides $\mathcal{F}\ell(\gamma)$ the structure of $U(\mathfrak{gl}_n)$-module, and respects the relations within $\mathcal{F}\ell(\gamma)$ (cf. \cite{Ful:1}),  i.e. the module in question is actually a \emph{module algebra} for the Hopf algebra $U(\mathfrak{gl}_n)$. 

Now seems a good time to mention another important fact about the flag algebra, Hodge's basis theorem \cite{Hod:1,HodPed:1}.

\begin{mytheo}\label{thm:comm-basis-thm} The algebra $\mathcal{F}\ell(\gamma)$ has $\C$-basis given by the monomials $f_T := f_{I_1}\cdots f_{I_p}$ where $|I_1|\geq |I_2|\geq \cdots \geq |I_p|\in\Vert\gamma\Vert$ and $(I_1,\ldots, I_p)$ produce a semi-standard Young tableau $T$ when filling out (in the obvious manner) the $p$ columns of a Young diagram (of appropriate shape).
\end{mytheo}

Among the tableaux mentioned in the theorem, we focus on those coming from $(I_1, I_2, \ldots, I_p) = ([d_1], [d_2], \ldots, [d_p])$ for integers $d_1\geq d_2 \geq \cdots \geq d_p \in \Vert\gamma\Vert$, i.e. the first row of $T$ is full of $1$'s, the second row, $2$'s, \ldots, the $d_1$th row, $d_1$'s.
Notice that for all $i$, $E_{ii}\cdot f_T = \lambda_i f_T$ for some $\lambda_i$ (precisely, the $i$-content $c_i(T)$ of $T$, i.e. the length of the $i$th row of $T$). Moreover, $E_{ij}\cdot f_T = 0$ for $i<j$. One ultimately deduces the
\begin{mytheo}\label{thm:comm-module-model} The flag algebra $\mathcal{F}\ell(\gamma)$ is the direct sum (with multiplicity one) of highest-weight modules $M^{\lambda}$, as $\lambda$ runs over all partitions with at most $n$ parts and column-lengths $d_1,\ldots, d_p \in \Vert\gamma\Vert$.
\end{mytheo}


\section{The Yangian Setting}\label{sec:yangian}
We recall the definitions of the Yangian and its determinant. For more details about the origin and construction of the Yangian, including the useful $R$-matrix formalism, cf. \cite{Mol:1}.

\begin{mydefi}\label{def:yangian} The \emph{Yangian} for $\mathfrak{gl}_n$ is the complex, associative, unital algebra $Y_n$ with countably many generators
$t^{(1)}_{ij}$, $t^{(2)}_{ij} , \ldots$ where $1 \leq i,j \leq n$,
and defining relations
\begin{equation}\label{eq:def-rels}
[t_{ij}^{(r+1)},t_{kl}^{(s)}]-
[t_{ij}^{(r)},t_{kl}^{(s+1)}]=
t_{kj}^{(r)}t_{il}^{(s)}-t_{kj}^{(s)}t_{il}^{(r)},
\end{equation}

where $r,s=0,1,2,\ldots\;$ and $t_{ij}^{(0)}:=\delta_{ij}\cdot 1$.
\end{mydefi}
Collecting the generators $t^{(r)}_{ij}\,(r=0,1,\ldots)$ together in the generating series
\begin{equation}\label{eq:gen-series}
t_{ij} (u) = \delta_{ij} + t^{(1)}_{ij} u^{-1} + t^{(2)}_{ij}u^{-2} +
\cdots \in Y_n[[ u^{-1} ]]\,,
\end{equation}
we may express the relations more compactly.
\begin{myprop}
The system of relations (\ref{eq:def-rels}) is equivalent to:
\begin{equation}\label{eq:def-rels-series}
[t_{ij} (u), t_{kl} (v) ] =  {1\over u-v} \left(t_{kj} (u) t_{il} (v) -
t_{kj} (v) t_{il} (u)\right) ,
\end{equation}
where $1\leq i,j,k,l \leq n$, calculations being carried out in $Y_n[u,v][[u^{-1},v^{-1}]]$.
\end{myprop}

We next explain how to view $Y_n$ as a deformation of $\C[\mathrm{GL}_n]$ instead of as a deformation of $U(\mathfrak{gl}_n)$. In spite of its flaws, this point-of-view manages to bear some fruit in subsequent sections.

Collecting the generating series together as a matrix of generators $T(u) = (t_{ij}(u))$, we are reminded of the coordinate algebra for $\mathrm{GL}_n$: $Y_n$ is $\C\langle T(u)\rangle$ modulo something or other, just as $\C[\mathrm{GL}_n]$ is $\C\langle X\rangle$ modulo something or other. What's more, $Y_n$ even has a determinant, like $\C[\mathrm{GL}_n]$ does. 

\subsection{Determinants}
Given $t_{ij}(u)$ as in (\ref{eq:gen-series}), define $t_{ij}(u+a)$, $a\in \mathbb Z$, to be the power series in $u^{-1}$ obtained by expanding the various factors $(u+a)^{-p}=\sum_{0\leq q}\binom{p}{q} a^q u^{-q-p}$ appearing below.
$$
t_{ij}(u+a) = \delta_{ij} + \sum_{1\leq p} t_{ij}^{(p)}(u+a)^{-p} = \delta_{ij} + \sum_{1\leq p} \bigg(\sum_{q+q'=p} \binom{q}{q'}a^{q'} \bigg) u^{-p}
$$
\begin{mydefi}[Yangian Determinant]\label{def:yangian-det} For all $I,J\in\binom{[n]}{d}$, the \emph{Yangian determinant} $\Det\!{}_{I,J}$ of  $T(u)$ is the power series in $u^{-1}$ given by the formula
\begin{eqnarray}
\nonumber \Det\!{}_{I,J} T(u) &=& t^I_J(u) = \sum_{\sigma\in \Sym_d} (-1)^{\ell(\sigma)} \times \\
\label{eq:yangian-det-rows} && t_{i_{\sigma(1)} j_1}(u) t_{i_{\sigma(2)} j_2}(u-1) \cdots t_{i_{\sigma(d)} j_d}(u-d+1)\,.
\end{eqnarray}
\end{mydefi}
Using the defining relations for $Y_n$, one discovers the 

\begin{myprop} For all $I,J$ as above, the Yangian determinant is also expressible as a sum of column permutations, namely
\begin{equation}\label{eq:yangian-det-cols}
t^I_J(u) = \sum_{\sigma\in \Sym_d} (-1)^{\ell(\sigma)} t_{i_1j_{\sigma(1)}}(u-d+1) \cdots t_{i_{d-1} j_{\sigma(d-1)}}(u-1) t_{i_dj_{\sigma(d)}}(u)\,.
\end{equation}
Moreover, the Yangian determinant is alternating in rows and columns, i.e., for all $\sigma\in\Sym_d$, for all $I, J\in[n]^d$, one has
\begin{equation}
\label{eq:ydet-alternating} t^{\sigma I}_{J} (u) = (-1)^{\ell(\sigma)} t^I_{J}(u) = t^I_{\sigma J}(u) \\
\end{equation}
using the right-hand side of either (\ref{eq:yangian-det-rows}) or (\ref{eq:yangian-det-cols}) to extend the definition of $\Det\!$ from sets to tuples.
\end{myprop}

For all $\alpha\in \C$, let $t^I_J(u+\alpha)$ denote the power series obtained by replacing each $t_{ij}(u-k)$ appearing on the right-hand side of (\ref{eq:yangian-det-rows}) by $t_{ij}(u-k+\alpha)$ before expanding. It is important to note that the mapping $Y_n\rightarrow Y_n$ represented by $T(u)\mapsto T(u+\alpha)\, (\forall \alpha\in\C)$ is an algebra automorphism \cite{Mol:1}. Thus, every $\Det\!$-minor identity appearing in the sequel may be rewritten in many ways by replacing any $t^I_J(u+b)$ appearing therein by $t^I_J(u+b+\alpha)$. 

\begin{myprop}[Laplace Expansion]\label{thm:yangian-Laplace-expansion} For all $d$-tuples $I, J\in[n]^d$, and all $1\leq r< d$ one has the \emph{cofactor expansion} relations
\begin{equation}\label{eq:ydet-row-expansion}
\sum_{\Lambda\subseteq[d],\,|\Lambda|=r} (-1)^{\Ell{\Lambda}{\,[d]\setminus\Lambda}} t^{I_{[r]}}_{J_{\Lambda}}(u-d+r)\cdot t^{I\setminus I_{[r]}}_{J\setminus J_{\Lambda}}(u) = t^I_J(u).
\end{equation}
\begin{equation}\label{eq:ydet-col-expansion}
\sum_{\Lambda\subseteq[d],\,|\Lambda|=r} (-1)^{\Ell{\Lambda}{\,[d]\setminus\Lambda}} t^{I_{\Lambda}}_{J_{[r]}}(u) \cdot t^{I\setminus I_{\Lambda}}_{J\setminus J_{[r]}}(u-d+r) = t^I_J(u).
\end{equation}
\end{myprop}
\begin{myprop}\label{thm:yangian-commuting-minors} For all $a\in\mathbb Z$, and for all subsequences $I', J'$ of $I$ and $J$, one has the \emph{commuting relation} $\left[t^I_J(u)\,,\, t^{I'}_{J'}(u+a)\right] = 0$. Equivalently, 
\begin{equation}\label{eq:yangian-comm}
\left[t^I_J(u)\,,\, t^{I'}_{J'}(v)\right] = 0.
\end{equation}
\end{myprop}

\subsection{ Coordinates}
Leaving for the moment the question of geometry (i.e. Yangian flags), let us follow Proposition \ref{thm:comm-flag-is-subalg} in an effort to define the algebraic counterpart (i.e. its ring of coordinate functions). 
Do we study the subalgebra of $Y_n[[u^{-1}]]$ generated by $\mathcal{G}_{\gamma}=\left\{t^{[d]}_I(u) : |I|=d, d\in\Vert\gamma\Vert\right\}$? If we want a $Y_n$ module (or comodule) structure, it is better to study the subalgebra of $Y_n$ generated by the coefficients of powers of $u^{-1}$ appearing in $\mathcal{G}_{\gamma}$. We are left with trying to find all of the relations among these ``coordinates,'' so we may give an abstract definition in terms of generators and relations as in Definition \ref{def:cplx-flag-algebra}. Also, we must describe the comodule and module structures. Toward the former goal, we have the following proposition, cf. \cite{Mol:2,Mol:1}.

\begin{myprop} For all tuples $A,I\in[n]^e$ and $B,J\in[n]^d$ with $e\geq d$, 
\begin{eqnarray}
\nonumber \big[t^A_I(u) \,,\,t^B_J(v) \big] &=& \sum_{p=1}^d \frac{(-1)^{p-1}p!}{(u-v-e+1)\cdots (u-v-e+p)} \times \\
\nonumber &&\Bigg(\sum_{\newatop{1\leq k_1<\cdots <k_p\leq n}{1\leq \ell_1<\cdots <\ell_p\leq n}} t^{a_1\cdots b_{\ell_1}\cdots b_{\ell_p}\cdots a_e}_{i_1\cdots i_e}(u) \cdot t^{b_1\cdots a_{k_1}\cdots a_{k_p}\cdots b_d}_{j_1\cdots j_d}(v) \\
\label{eq:yangian-C-rels}&& - t^{a_1\cdots a_e}_{i_1\cdots j_{\ell_1}\cdots j_{\ell_p}\cdots i_e}(v)\cdot t^{b_1\cdots b_d}_{j_1\cdots i_{k_1}\cdots i_{k_p}\cdots j_d}(u)\Bigg)
\end{eqnarray}
is a relation among the minors holding in $Y_n[u,v][[u^{-1},v^{-1}]]$. 
\end{myprop}

In the next section we find more, but for now notice that (\ref{eq:yangian-comm}) is a special case of (\ref{eq:yangian-C-rels}) (after the alternating property).


\def\StablyFinite{A ring is called \emph{Dedekind finite} if $\forall c,d\in R$, $ab=1$ implies $ba=1$. Such a ring is called \emph{stably finite} if this property continues to hold for the $d\times d$ matrices over $R$ with $d>1$. See \cite{Lam:2} for more details, and \cite{Coh:4, Mon:1} for some non-stably finite rings.}
\def\qdetInduction{A proof of the quantum-determinantal analog of this claim appears in \cite{Lau:1}} 
\section{Noncommutative Flags}\label{sec:NC-flags}
Given a composition $\gamma\models n$ and a skew field $D$, one may define the noncommutative flags $Fl(\gamma)$ as in Section \ref{sec:comm-flags} (the two choices for definition agree by the \emph{invariant basis number} property of skew fields, cf. \cite{Lam:1}). 
Clearly one cannot use the determinant to coordinatize the noncommutative flags, but the quasideterminant of Gelfand and Retakh (\cite{GelRet:1, GGRW:1}) offers an alternative (\cite{GelRet:3, Lau:4}). 
In this section, we describe the quasi-Pl\"ucker coordinates, the relations known to hold among them, and how they specialize in the Yangian setting. 

\subsection{Quasi-Pl\"ucker Coordinates}\label{sec:qplucker-coords}
We assume in this subsection that $A=(a_{ij})$ is a matrix of noncommuting 
indeterminants and that $D$ is the \emph{free skew field} $F\lskew A\rskew$ (cf, e.g., \cite{Coh:2, GelRet:3}). 
We do this in lieu of taking a \emph{generic} point in $D^{n^2}$ for arbitrary $D$, by which we mean ``every submatrix we wish to invert is indeed invertible.'' When ``specializing'' to arbitrary division rings, this means that the coordinates and equations in this section will make sense on a dense subset of the big Schubert cell in $Fl(\gamma)$.

\begin{mydefi} For each $i,j\in[n]$ we define the $(i,j)^{\mathrm{th}}$ quasideterminant $|A|_{ij}$ of $A$ by the formula
$$
|A|_{ij} = a_{ij} - \rho_{ij} \cdot (A^{ij})^{-1} \cdot \chi_{ji}\,,
$$
where $\rho_{ij} = A_{\{i\},[n]\setminus j}$ is row $i$ of $A$ with column $j$ deleted, and $\chi_{ji} = A_{[n]\setminus i, \{j\}}$ is column $j$ of $A$ with row $i$ deleted.
\end{mydefi}

A simple calculation shows
\begin{equation}\label{eq:qdet-is-inv}
|A|_{ij}^{-1} = (A^{-1})_{ji},
\end{equation}
which continues to hold over \textit{less-free} skew fields $D$ provided both sides are defined. Extending the definition to submatrices $A_{I,J}$ of $A$ in the obvious fashion, we have the important

\begin{mytheo}[Column Homological Relations, \cite{GelRet:1}]\label{thm:colhomol-rels} For any $L,M\subseteq [n]$, $i,j\in[n]$ with $|L|=|M|+1$ and $i,j\not\in M$ we have $(\forall s\not\in L, \forall t\in L)$:
\begin{equation}\label{eq:colhomol-rels}
|A_{sL,ijM}|_{si}\cdot |A_{L,iM}|_{ti}^{-1} = - |A_{sL,ijM}|_{sj}\cdot |A_{L,jM}|_{tj}^{-1}\,.
\end{equation}
\end{mytheo}

As an immediate corollary, one sees that the left ratio $|A_{L, iM}|_{ti}^{-1} |A_{L, jM}|_{tj}$
is independent of the choice of $t\in L$. This allows us to make the

\begin{mydefi}[Left/Column Coordinates] For $n$, $A$, and $M$ as above, the \emph{quasi-Pl\"ucker coordinate} of \emph{order} $|M|+1$ associated to $(i,j,M)$ is given by
$$
p_{ij}^M=p_{ij}^M(A) = |A_{[m], iM}|_{si}^{-1} |A_{[m], jM}|_{sj} \quad(\hbox{any }s\in[m]).
$$
\end{mydefi}
\begin{myrema} Regarding the ``left'' and ``column'' tags appearing above: there is a symmetric theory involving various row sets of $A$ and considering right ratios up to equivalence by a right action of $\mathrm{P}_{\gamma}^{+}$ (the block upper-triangular parabolic subgroup). We will not need this in the present paper.
\end{myrema}

\begin{mytheo}[\cite{GelRet:3}] \label{thm:quasi-gln-invariant} Fix $i,j,m,n,M$, and $A$ as above. Put $B=A_{[m],[n]}$. For any $g\in \mathrm{GL}_m(D)$, 
$$p_{ij}^M(g\cdot B) = p_{ij}^M(B)\,.$$
\end{mytheo}

We apply these constructions to our problem of coordinatizing flags by taking $m\in\Vert\gamma\Vert$, and viewing $A$ as some $A(\Phi) \in Fl(\gamma)$. After Theorem \ref{thm:quasi-gln-invariant}, we learn that quasi-Pl\"ucker coordinates are not projective invariants of $A$, but true invariants). Still, the set $\{p_{ij}^M \mid |M|+1\in\Vert\gamma\Vert \}$ describes $Fl(\gamma)$ in the following sense: (i) no greater collection of quasi-Pl\"ucker coordinates is invariant under $\mathrm{P}_{\gamma}$; (ii) if $f$ is a function on $A$ which is $\mathrm{P}_{\gamma}$ invariant, then $f$ is a rational function on the collection $p_{ij}^M(A)$.

Working toward a statement analogous to Theorem \ref{thm:comm-coords-characterize}, we start with noncommutative analogs of the alternating and Young symmetry relations:

\begin{mytheo}\label{thm:quasi-alternating} Let $A=(a_{ij})$ be an $n\times n$ matrix of formal, noncommuting variables. The following identities hold in $F\lskew A\rskew$:
\begin{itemize}
\item Fix $M\in[n]^d\, (d<n)$ with distinct entries. If $i,j\in[n]$ with $i\not\in M$, then putting $B=A_{[d+1],(i| j|M)}$, we have
\begin{center}$p_{ij}^M(B)$ does not depend on the ordering of $M$. 
\end{center}

\item Fix $M\in\binom{[n]}{d}\, (d<n-2)$. If $i,j,k\in[n]\setminus M$, then we have
$$p_{ij}^{k\cup M}\, p_{jk}^{i\cup M}\,p_{ki}^{j\cup M} = -1\,.$$

\item Fix $M\in\binom{[n]}{d}\, (d<n)$. If $i,j\in[n]$ with $i\not\in M$, then putting $B=A_{[d+1],(j|i\cup M)}$, we have
$$
p_{ij}^M(B)=\left\{\begin{array}{ll}
0 & \hbox{if }j\in M\\
1 & \hbox{if }j=i\end{array}\right. .$$

\item Fix $M \in\binom{[n]}{d}\, (d<n-1)$. If $i,j,k\in[n]$ with $i,j\not\in M$, then we have
$$ p_{ij}^M p_{jk}^M = p_{ik}^M\,.
$$
\end{itemize}
\end{mytheo}

\begin{mytheo}[Quasi-Pl\"ucker Relations]\label{thm:quasi-YS} Let $A$ be an $n\times n$
matrix of formal, noncommuting variables. Fix $L,M\in[n]$ with $s=|L|\geq|M|+1=t$ and $s,t\in\Vert\gamma\Vert$. Fix $i\in[n]\setminus M$. The following identities hold in $F\lskew A\rskew$
\begin{equation}\label{eq:quasi-YS-left}
\sum_{j\in L}p_{ij}^{M}\cdot p_{ji}^{L\setminus j} = 1\,.
\end{equation}
 \end{mytheo}
This was observed for the case $|M|+1=|L|$ in \cite{GelRet:3}. We abbreviate these relations as $({\mathcal P}_{i,M,L})$. See \cite{Lau:1} for a complete proof in the symmetric case involving row coordinates.

Unfortunately, it is not known if the previous two Theorems ``exhaust'' the fundamental relations holding among the quasi-Pl\"ucker coordinates, i.e. a noncommutative version of Theorem \ref{thm:comm-coords-characterize} remains elusive. Still, there does exist the following very compelling prelude:

\begin{mytheo}\label{thm:g-invt-fcns}
Let $A=(a_{ij})$ be an $n\times n$ matrix with formal, noncommuting entries and suppose $f=f(a_{ij})$ is a rational function over the free skew-field $D=F\lskew A\rskew$. If $f(gA)=f(A)$ for all $g\in\mathrm{P}_{\gamma}(D)$, then $f$ is a rational function in the quasi-Pl\"ucker coordinates $\{p_{ij}^M(A) \,:\, |M|+1\in\Vert\gamma\Vert \}$.
\end{mytheo}

A Grassmannian version of this theorem appears in \cite{GelRet:3}. The proof is a consequence of noncommutative Gaussian Elimination and a simple application of the noncommutative Sylvester's Identity (\cite{GGRW:1}) and induction. We illustrate the theorem with a $3\times 3$ example, $\gamma=(2,1)$.

\begin{proof}[Sketch of Proof]
Using only elements of $\mathrm{P}_{\gamma}$, we may transform $A$ into
$$
\left[\begin{array}{ccc} 1 & a_{11}^{-1}a_{12} & a_{11}^{-1}a_{13} \\
0 & |A_{\{1,2\},\{1,2\}}|_{22} & |A_{\{1,2\},\{1,3\}}|_{23} \\
0 & |A_{\{1,3\},\{1,2\}}|_{32} & |A_{\{1,3\},\{1,3\}}|_{33}
\end{array}\right],
$$
and into
$$
\left[\begin{array}{ccc} 1 & a_{11}^{-1}a_{12} & a_{11}^{-1}a_{13} \\
0 & 1 & |A_{\{1,2\},\{1,2\}}|_{22}^{-1} |A_{\{1,2\},\{1,3\}}|_{23} \\
0 & 0 & |A_{\{1,2,3\},\{1,2,3\}}|_{33}
\end{array}\right].
$$
Continuing Gaussian Elimination via elements of $\mathrm{P}_{\gamma}$, we reach the matrix
\begin{equation}\label{eq:LDU}
\left[\begin{array}{ccc} 1 & 0 & p_{13}^{\emptyset}-p_{12}^{\emptyset}p_{23}^{1} \\
0 & 1 & p_{23}^{1} \\
0 & 0 & 1
\end{array}\right].
\end{equation}
Consequently, $f$ is a rational function in the Pl\"ucker coordinates $p_{ij}^M$ of $A$. However, not all $M$ appearing satisfy the hypotheses of the theorem; e.g. the symbol $p_{13}^{\emptyset}$ it is of order $1$, while the allowable orders are $\Vert\gamma\Vert=\{2\}$. We have a little more work to do. From Theorems \ref{thm:quasi-alternating} and \ref{thm:quasi-YS}, we see that 
\begin{eqnarray*}
p_{13}^{\emptyset}-p_{12}^{\emptyset}p_{23}^{1} &=& (p_{13}^{\emptyset}p_{31}^{2} - p_{12}^{\emptyset}p_{23}^{1}p_{31}^{2})p_{13}^{2} \\
&=& (p_{13}^{\emptyset}p_{31}^{2} + p_{12}^{\emptyset}p_{21}^{3})p_{13}^{2}\\
&=& p_{13}^{2}\,,
\end{eqnarray*}
so we are left with the reduced form of $A$ looking like
\begin{equation}\label{eq:LDU-std}
\left[\begin{array}{ccc} 1 & 0 & p_{13}^{2} \\
0 & 1 & p_{23}^{1} \\
0 & 0 & 1
\end{array}\right].
\end{equation}
In short, if $\gamma=(\gamma_1,\ldots,\gamma_r)$, then rows $|\gamma_{[i-1]}|+1$ through $|\gamma_{[i]}|$ of the reduced form of $A$ will consist of a $\gamma_i \times |\gamma_{[i-1]}|$ block of zeros beside an identity matrix (of order $\gamma_i$) beside a collection of left quasi-Pl\"ucker coordinates of order $|\gamma_{[i]}|$, $1\leq i\leq r$.
\end{proof}

\subsection{$T$-generic Flags}\label{sec:Tgeneric-flags}
Consider an algebra $\mathcal A(n)$ on $n^2$ generators $t_{ij}$ over a field $F$---ignoring the relations for now. Suppose $\mathcal A(n)$ may be embedded in $D$, and put all the generators together in a matrix $T$. 
We view the $t_{ij}$ as coordinate functions and their relations as characterizing some set $X$ inside $D^{n^2}$. Let us call $X$ the set of \emph{$T$-generic matrices over $D$ for $\mathcal A(n)$}. By the \emph{$T$-generic flags over $D$ for $\mathcal A(n)$} we mean those cosets in $\mathrm{P}_{\gamma} \,\backslash\, \mathrm{GL}_n(D)$ having a representative in $X$. If $T$ is invertible, this set is evidently nonempty. If its submatrices are also invertible, then all of the identities displayed above carry over to the $T$-generic setting. 

One might then study the quasi-Pl\"ucker coordinates of $T$ toward describing a flag algebra for $\mathcal{A}(n)$. An important example of a setting $\mathcal{A}(n)$ where the above is possible is the quantum group $\mathrm{GL}_q(n)$ and its flag algebra \cite{TafTow:1, Lau:1}. 
The Yangian does not quite fit into this rubrick, but it comes close, and close enough to help us define $\mathcal{F}\ell(\gamma)$ for $Y_n$. 

In \cite{TafTow:1}, Taft and Towber find three types of relations among the quantum minors and go on to show that these three are sufficient to give the ``correct'' quantum generalization to the algebra $\mathcal{F}\ell(\gamma)$ of Section \ref{sec:comm-flags}. The first two are quantum versions of the alternating and Young symmetry relations outlined in (\ref{eq:cplx-alternating}) and (\ref{eq:cplx-YS}). The third type---which we shall call \emph{monomial straightening} relations---is a replacement for the commuting property $\big[f_I\,,\, f_J\big] = 0$ holding among coordinate functions of a (commutative) algebraic variety. Section \ref{sec:yangian-newrels} is dedicated to finding Yangian versions of the latter two types (the alternating relations being already given in (\ref{eq:ydet-alternating})). We first verify the hypotheses in the first paragraph above.

\subsection{Application to $Y_n$}
We want to show that $Y_n[[u^{-1}]]$ may be embedded in a division ring $D$. 
With a little effort one can rewrite the relations (\ref{eq:def-rels}) as follows
\begin{equation}\label{eq:yangian-filtered}
t^{(r)}_{ij} t^{(s)}_{kl}  =  t^{(s)}_{kl} t^{(r)}_{ij} + \sum_{a=1}^{\min(r,s)} \left( t^{(a-1)}_{kj} t^{(r+s-a)}_{il} - t^{(r+s-a)}_{kj} t^{(a-1)}_{il}\right).
\end{equation}

Call a word $t_{i_1j_1}^{(r_1)}t_{i_2j_2}^{(r_2)}\cdots t_{i_pj_p}^{(r_p)}$ in the generators a monomial of \emph{degree} $p$ and \emph{weight} $r_1 + r_2 + \cdots + r_p$. Then (\ref{eq:yangian-filtered}) says that $Y_n$ is a filtered algebra by weight. Moreover, (\ref{eq:yangian-filtered}) reveals that the associated graded algebra $\mathrm{gr}\hbox{-}{Y}_n$ is the commutative $\C$-algebra freely generated by the $\mathrm{gr}\hbox{-}Y_n$ image of the set $\left\{ t_{ij}^{(r)} \mid r\geq 1; 1\leq i,j\leq n\right\}$. In particular,  $\mathrm{gr}\hbox{-}{Y}_n$ is a (right) Ore domain, i.e., for all $x,y\in \mathrm{gr}\hbox{-}{Y}_n$, there exist $a,b\in \mathrm{gr}\hbox{-}{Y}_n$ such that $xa=yb$.

\begin{mytheo}[Cohn, \cite{Coh:3}] If $R$ is a filtered ring, and $\mathrm{gr}\hbox{-}R$ is a right (or left) Ore domain, then $R$ is embeddable in a skew field.
\end{mytheo}

Let $D_0$ denote the skew field for $Y_n$ provided by Cohn's theorem. It is easy to see that if $R$ is a skew field, then $R[[x]]$ may be embedded in a skew field as well, namely the Laurent series in $x$. Let $D$ be the division ring built in this way from $D_0[[u^{-1}]]$. 

Next, we must show that $T(u)$ and its submatrices are invertible over $D$. For this, we turn to the minors $t^I_J(u)$ of Section \ref{sec:yangian}. These minors are invertible in $D$---otherwise, they are zero, which is clearly not the case by virtue of their nonvanishing in $\mathrm{gr}\hbox{-}Y_n[[u^{-1}]]$. After Proposition \ref{thm:yangian-Laplace-expansion}, it is easy to see that $S_{I,J}(u) \cdot T_{I,J}(u) = 1$ in $D$, where for any $I,J\in\binom{[n]}{d}$, $S_{I,J}(u) = (s_{I,J}(u)_{kl})$ is the $d\times d$ matrix given by
\begin{equation}\label{eq:ydet-is-inv}
s_{I,J}(u)_{kl} = \frac{(-1)^{k+l}}{t^I_J(u+d-1)} \cdot t^{i_1\cdots \widehat{i_l}\cdots i_d}_{j_1\cdots \widehat{j_k}\cdots j_d} (u+d-1),
\end{equation}
and the factors on the right commute by (\ref{eq:yangian-comm}). There are rings $R$ such that an equality $ST=1$ concerning two square matrices over $R$ does not imply $TS=1$, but division rings are not among them.\footnote{\StablyFinite} Deduce, as desired, that $T_{I,J}(u)$ is invertible for all $I,J\in\binom{[n]}{d}$. 

We conclude with two equations that will be useful in the sequel: (i) using (\ref{eq:qdet-is-inv}) and (\ref{eq:ydet-is-inv}), we may write 
\begin{equation}\label{eq:ydet-qdet-factors}
p_{ab}^K(T(u)) = t^{[d]}_{a|K} (u+d-1)^{-1} \cdot t^{[d]}_{b|K}(u+d-1)
\end{equation} 
for any $K\in[n]^{d-1}$ with $a\in [n]\setminus K$; (ii) the homological relations (\ref{eq:colhomol-rels}) imply
\begin{equation}\label{eq:ydet-weak-qcomm}
t^L_{a|K}(u) \cdot t^L_{b|K}(u+1) = t^L_{b|K}(u) \cdot t^L_{a|K}(u+1)\,.
\end{equation}
View (\ref{eq:ydet-weak-qcomm}) as a weak version of the well-known $q$-commuting relations holding for the quantum determinant, cf. \cite{KroLec:1, Lau:1}.

\subsection{New Relations}\label{sec:yangian-newrels}
\begin{mylemm} For all $1\leq d\leq e < n$, and $I, J\subseteq [n]$ with $|I|=d-1, |J|=e+1$, the Yangian minors of the matrix $T(u)$ satisfy
\begin{equation}\label{eq:yangian-YS-rels}
0 = \sum_{\lambda\in J} (-1)^{\Ell{\lambda}{J\setminus \lambda}} t^{[d]}_{I|\lambda}(u+d) \cdot t^{[e]}_{J\setminus \lambda}(u+e+1).
\end{equation}
\end{mylemm}
\begin{proof} A straightforward exercise in clearing denominators, starting from (\ref{eq:quasi-YS-left}) and using (\ref{eq:ydet-qdet-factors}) and (\ref{eq:ydet-weak-qcomm}). 
\end{proof}

This looks like the Young symmetry relations of (\ref{eq:cplx-YS}) except we only know the result for $|\Lambda|=1$. Two generalizations may be proposed, and both are true. 

\begin{myprop}[Young Symmetry]\label{thm:yangian-YS} For all $0<p\in\mathbb N$,  and all $I\in[n]^{d-p}$ and $J\in\binom{[n]}{e+p}$ with $d\leq e$, the Yangian minors of $T(u)$ satisfy
\begin{equation}\label{eq:yangian-YS-rels-I}
0 = \sum_{{\Lambda\in J},\,{|\Lambda|=p}} (-1)^{\Ell{\Lambda}{J\setminus \Lambda}} t^{[d]}_{I|\Lambda}(u+d) \cdot t^{[e]}_{J\setminus \Lambda}(u+e+1) 
\end{equation}
and
\begin{equation}\label{eq:yangian-YS-rels-II}
0 = \sum_{{\Lambda\subseteq J},\,{|\Lambda|=p}} (-1)^{\Ell{\Lambda}{J\setminus \Lambda}} t^{[d]}_{I|\Lambda}(u+d) \cdot t^{[e]}_{J\setminus \Lambda}(u+e+p)\,.
\end{equation}
\end{myprop}

\begin{proof}[Proof of (\ref{eq:yangian-YS-rels-I})] 
Let ${Y_{I, J}}_{(p)}$ represent the right-hand side of (\ref{eq:yangian-YS-rels-I}). We claim that\footnote{\qdetInduction}
$$ 
{Y_{I, J}}_{(p)} = \frac{1}{\,p\,}\sum_{\lambda\in J} (-1)^{\Ell{\lambda}{J\setminus \lambda}}\cdot {Y_{I|\lambda,\, J}}_{(p-1)},
$$
which verifies Equation (\ref{eq:yangian-YS-rels-I}) after the lemma and induction on $p$.
\end{proof}

\begin{proof}[Proof of (\ref{eq:yangian-YS-rels-II})]
We demonstrate (\ref{eq:yangian-YS-rels-II}) by combining Laplace expansions of $\Det\!$. According to (\ref{eq:ydet-col-expansion}), we have
\begin{equation*}
t^{[d]}_{I |\Lambda}(u+d) = \sum_{K\subseteq [d],\, |K|=|\Lambda|} (-1)^{\Ell{[d]\setminus K}{K}} t^{[d]\setminus K}_{I}(u+d) t^{K}_{\Lambda}(u+d-(d-p))\,,
\end{equation*}
while (\ref{eq:ydet-row-expansion}) tells us
\begin{equation*}
\sum_{\Lambda} (-1)^{\Ell{\Lambda}{J\setminus \Lambda}} t^{K}_{\Lambda}(u+p)) t^{[e]}_{J\setminus\Lambda}(u+e+p) = t^{K\, |\,[e]}_{J}(u+e+p)\,.
\end{equation*}
Since $[d]\subseteq [e]$, the tuple $(K\,|\,[e])$ has repeated indices. Conclude that the right-hand side of (\ref{eq:yangian-YS-rels-II}) is a sum over $K$ with summand identically zero.
\end{proof}

\begin{myrema} To prove the lemma preceding Proposition \ref{thm:yangian-YS}, one starts from a equation of the form $1=\sum_{\lambda} \left({t_{M|i}}^{-1}t_{M|\lambda}\right)\left({t_{L\setminus \lambda|\lambda}}^{-1} t_{L\setminus j|i}\right)$ for carefully chosen $(i,M,L)$ and ``clears denominators'' in two different directions. What do we learn if we move both inverted minors to the other side in the same direction? 
\end{myrema}

Apply $(\mathcal P_{i,M,L})$ to $T(u+1)$, putting $M=\emptyset$, and deduce:
\begin{eqnarray*}
1 &=& \sum_{\lambda\in L} t^{[1]}_{i}({\scriptstyle u+1})^{-1} t^{[1]}_{\lambda}({\scriptstyle u+1}) t^{[e]}_{L\setminus \lambda|\lambda}({\scriptstyle u+e})^{-1} t^{[e]}_{L\setminus \lambda|i}({\scriptstyle u+e}) \\
&=& \sum_{\lambda\in L} (-1)^{\Ell{L\setminus \lambda}{\lambda}} t^{[1]}_{i}({\scriptstyle u+1})^{-1} t^{[e]}_{L}({\scriptstyle u+e})^{-1} t^{[1]}_{\lambda}({\scriptstyle u+1}) t^{[e]}_{L\setminus \lambda|i}({\scriptstyle u+e})\,, \\
\end{eqnarray*}
or
$$
t^{[e]}_{L}({u+e}) t^{[1]}_{i}({u+1}) = \sum_{\lambda\in L} (-1)^{\Ell{L\setminus \lambda}{\lambda}} t^{[1]}_{\lambda}({u+1}) t^{[e]}_{L\setminus \lambda|i}({u+e})\,.
$$
This equation may be understood as a straightening relation: \emph{among all monomials involving minors of different orders, we prefer those whose minors are arranged in ascending order.} The case $M=\emptyset$ suggests a general phenomenon.

\begin{myprop}[Monomial Straightening] For all $I,J\subseteq[n]$ with $|I|=d\leq e=|J|$, 
\begin{equation}\label{eq:yangian-MS-rels}
t^{[e]}_{J}({u+e}) t^{[d]}_{I}({u+d}) = \sum_{\newatop{\Lambda\subseteq J}{|\Lambda|=e-d}} (-1)^{\Ell{\Lambda}{J\setminus \Lambda}} t^{[d]}_{J\setminus \Lambda}({u+d}) t^{[e]}_{\Lambda|M}({u+e})\,.
\end{equation}
\end{myprop}

\begin{proof} A calculation analogous to the demonstration of (\ref{eq:yangian-YS-rels-II}), mixing the two Laplace expansions of $\Det\!$. 
\end{proof}

\section{Main Results}\label{sec:main-results}
We propose to study the following class of algebras as \emph{the Yangian flag algebras.}
\begin{mydefi} Given any composition $\gamma\models n$, let $\mathcal{F}\ell_{Y_n}\!(\gamma)$ be the $\C$-algebra with generators $\left\{ f_I^{(r)} \mid |I|\in \Vert\gamma\Vert,\, r=0,1,\ldots \right\}$ and the \emph{alternating}, \emph{commuting}, \emph{Young symmetry}, and \emph{monomial straightening} relations given below:
\begin{description}
\item[$(\mathcal A_J)$] $\forall J\in[n]^d,\, \forall\sigma\in\Sym_d$
\begin{equation}\label{eq:yang-flag-A-rels} f_{\sigma J} (u) = (-1)^{\ell(\sigma)} f_{J}(u) \,. 
\end{equation}

\item[$(\mathcal C_{I,J})$] $\forall I\in[n]^d$ and $J\in[n]^e$ with $1\leq d\leq e$
\begin{eqnarray}
\label{eq:yang-flag-C-rels} \big[f_J(u) \,,\,f_I(v) \big] &=& \sum_{p=1}^d \frac{(-1)^{p-1}p!}{(u-v-e+1)\cdots (u-v-e+p)} \times \\
\nonumber &&\Bigg( \binom{d}{p} f_J(u) f_I(v) \\
\nonumber && - \sum_{K,L\in\binom{[n]}{p}}  f_{i_1\cdots j_{\ell_1}\cdots j_{\ell_p}\cdots i_e}(v)\cdot f_{j_1\cdots i_{k_1}\cdots i_{k_p}\cdots j_d}(u)\Bigg).
\end{eqnarray}

\item[$(\mathcal Y_{I,J})$] $\forall 1\leq p\leq d\leq e,\,\forall I\in[n]^{d-p},\,\forall J\in\binom{n}{e+p}$
\begin{equation}\label{eq:yang-flag-YS-rels}
0 = \sum_{{\Lambda\subseteq J},\,{|\Lambda|=p}} (-1)^{\Ell{\Lambda}{J\setminus \Lambda}} f_{I|\Lambda}(u+d) \cdot f_{J\setminus \Lambda}(u+e+p)\,.
\end{equation}

\item[$(\mathcal M_{I,J})$] $\forall d\leq e,\, \forall I\in\binom{n}{d},\,J\in\binom{n}{e}$
\begin{equation}\label{eq:yang-flag-MS-rels}
f_{J}({u+e}) f_{I}({u+d}) = \sum_{\newatop{\Lambda\subseteq J}{|\Lambda|=e-d}} (-1)^{-\Ell{\Lambda}{J\setminus \Lambda}} f_{J\setminus \Lambda}({u+d}) f_{\Lambda|I}({u+e})\,.
\end{equation}
\end{description}
We understand the above equations as giving relations among the generators of $\mathcal{F}\ell_{Y_n}\!(\gamma)$ by introducing the power series $f_I(u):=f_I^{(1)}u^{-1}+f_I^{(2)}u^{-2}+\cdots$ and comparing coefficients of the different powers of $u^{-r}v^{-s}$ appearing on each side.
\end{mydefi}


\subsection{Comodules}
In addition to having a matrix of generators and a determinant, $Y_n$ shares another important feature with $\C[\mathrm{GL}_n]$, cf. \cite{Mol:1}.
\begin{mytheo} The Yangian $Y_n$ is a bialgebra with structure maps given by
$$
\Delta : t_{ij}(u) \mapsto \sum_{1\leq k\leq n} t_{ik}(u) \otimes t_{kj}(u) \qquad \varepsilon : T(u) \mapsto 1.
$$
\end{mytheo}
These expressions are to be understood as maps by sending, e.g., $t_{ij}^{(r)}$ to the coefficient of $u^{-r}$ appearing in the expansion of $t_{ik}(u) \otimes t_{kj}(u)$.

\begin{mycoro} The Yangian minors diagonalize according to the formula
\begin{equation}\label{eq:yangian-minor-Delta}
\Delta t^I_J(u) = \sum_{K\subseteq[n], |K|=|I|} t^I_K(u) \otimes t^K_J(u).
\end{equation}
\end{mycoro}

We are now ready to state our first main result.

\begin{mytheo} The flag algebra $\mathcal{F}\ell_{Y_n}\!(\gamma)$ is a right $Y_n$-comodule algebra with structure map given by $\rho(f_I(u) ) = \sum_{J\subset[n], |J|=|I|} f_J(u) \otimes t^J_I(u)$.
\end{mytheo}

Toward a proof, we begin by focusing on a model for the degree one, weight $d$ constituents of $\mathcal{F}\ell_{Y_n}\!(\gamma)$. We work with the full flag $\mathcal{F}\ell_{Y_n}\!(n)$ for simplicity, letting the reader supply the necessary changes for arbitrary $\gamma$.

\begin{mydefi} For all $1\leq d\leq n$, let $C_n(d)$ be the vector space spanned by $\{\tilde f_I^{(r)} \mid I\in{[n]}^{d},\, r=0,1, \ldots \}$ modulo the alternating relations $(\mathcal A_I)$:
\begin{equation*}
\tilde f_{\sigma I}(u) = (-1)^{\ell(\sigma)} \tilde f_I(u) \qquad (\forall I\in[n]^d,\,\forall \sigma\in\Sym_d),
\end{equation*}
arranging the generators in a power series as usual.
\end{mydefi}

\begin{mylemm} $C_n(d)$ is a right $Y_n$-comodule with structure map given by
$$
\rho\left(\tilde f_I(u)\right) = \sum_{J\in\binom{[n]}{d}} \tilde f_J(u) \otimes t^J_I(u)
$$
for all $I\in[n]^d$.
\end{mylemm}
\begin{proof} After (\ref{eq:yangian-minor-Delta}), we need only check that $\rho$ respects the relations. But this is evident after (\ref{eq:ydet-alternating}).
\end{proof}

Conclude that $C_n := \bigoplus_{1\leq d<n} C_n(d)$ is also a right $Y_n$-comodule. 
Extend $\rho$ to tensor products $C_n(d)\otimes C_n(e)$ in the usual way:
$$
\rho^{\otimes 2}\left(\tilde f_I(u)\otimes \tilde f_J(v)\right) = \sum_{K,L} \tilde f_K(u) \otimes \tilde f_L(v)\otimes 
t^K_I(u) t^L_J(v)\,,
$$
to deduce that the tensor algebra $T(C_n)$ is a right $Y_n$-comodule algebra. 

\begin{myremas} \textit{1.} Note that $\mathcal{F}\ell_{Y_n}\!(n)$ is $T(C_n)$ modulo the three sets of relations $(\mathcal Y_{I,J})$, $(\mathcal M_{I,J})$, and $(\mathcal C_{I,J})$. If we can show that $T(C_n)$ modulo each of these is again a comodule algebra, we will have proven that $\mathcal{F}\ell_{Y_n}\!(n)$ is as well.

\noindent\textit{2.} For any scalar $\alpha\in\C$, let $C_n({\alpha};d,e)\subseteq C_n(d)\otimes C_n(e)$ be the span of the coefficients of the powers of $v^{-r}\, (r\geq 1)$ appearing in $\left\{\tilde f_I(v+\alpha) \otimes \tilde f_J(v) \mid I\in\binom{[n]}{d}, J\in\binom{[n]}{e} \right\}$. The preceding equation reveals that $C_n({\alpha};d,e)$ is a subcomodule of $C_n(d)\otimes C_n(e)$.
\end{myremas}

\begin{myprop} The map $\sigma: C_n(e+p) \rightarrow C_n(p)\otimes C_n(e)$ given by
$\tilde f_J(u) \mapsto \sum_{\Lambda\subseteq J,\,|\Lambda|=p} (-1)^{\Ell{\Lambda}{J\setminus \Lambda}} \tilde f_\Lambda(u-e)\otimes \tilde f_{J\setminus \Lambda}(u)$ is a comodule map.
\end{myprop}

\begin{proof} The map is evidently only well-defined on $\tilde f_J(u),\,J\in\binom{[n]}{e}$; for instance $0 = \tilde f_{22}(u) \mapsto 2\tilde f_2(u-2)\otimes \tilde f_2(u)$. With this restriction in mind, we compare the action of ($\star$) $\rho^{\otimes 2}\circ\sigma$ and ($\star\star$) $(\sigma\otimes 1)\circ \rho$ on $\tilde f_J(u)$.
\begin{eqnarray*}
(\star)\,\,\, \tilde f_J({u}) &=& \sum_{\Lambda\subseteq J} (-1)^{\Ell{\Lambda}{J\setminus \Lambda}} \sum_{M,K_0} \tilde f_M({\scriptstyle u-e})\otimes \tilde f_{K_0}({\scriptstyle u}) \otimes t^M_{\Lambda}({\scriptstyle u-e}) t^{K_0}_{J\setminus \Lambda}({\scriptstyle u}) \\
&=& \sum_{M,K_0} \tilde f_M({\scriptstyle u-e})\otimes \tilde f_{K_0}({\scriptstyle u}) \otimes\left(\sum_{\Lambda\subseteq J} (-1)^{\Ell{\Lambda}{J\setminus \Lambda}}  t^M_{\Lambda}({\scriptstyle u-e}) t^{K_0}_{J\setminus \Lambda}({\scriptstyle u}) \right)\\
&=& \sum_{M,K_0} \tilde f_M({\scriptstyle u-e})\otimes \tilde f_{K_0}({\scriptstyle u}) \otimes t^{M|K_0}_{J}({\scriptstyle u}),
\end{eqnarray*}
using the row Laplace expansion (\ref{eq:ydet-row-expansion}). On the other hand,
\begin{eqnarray*}
(\star\star)\,\,\, \tilde f_J(u) &=& \sum_{K}\sum_{\newatop{M\subseteq K}{K_0=K\setminus M}} (-1)^{\Ell{M}{K_0}} \tilde f_M({\scriptstyle u-e}) \otimes \tilde f_{K_0} ({\scriptstyle u}) \otimes t^{K}_{J}({\scriptstyle u}) \rule[0pt]{8em}{0pt}\\
&=& \sum_{M,K_0} \tilde f_M({\scriptstyle u-e})\otimes \tilde f_{K_0}({\scriptstyle u}) \otimes t^{M|K_0}_{J}({\scriptstyle u}).
\end{eqnarray*}
\end{proof}

\begin{myprop} The map $\mu: C_n({d-p};d-p,p) \rightarrow C_n(d)$ given by
$\tilde f_I(u+d-p) \otimes \tilde f_\Lambda(u) \mapsto \tilde f_{I|\Lambda}(u+d-p)$ is a comodule map.
\end{myprop}
\begin{proof} Analogous to the preceding proof. This time the column Laplace expansion is used. Note also that this map is well-defined for $\tilde f_I(u),\,I\in[n]^{d}$.
\end{proof}

Recall that the composition of comodule maps is again a comodule map, and that the image of a comodule map is another comodule. This allows us to conclude that the tensor algebra $T(C_n)$ modulo the Young symmetry relations $(\mathcal Y_{I,J})$ is a comodule algebra. For if we apply $(\mu\otimes 1)\circ (1\otimes \sigma)$ to $\tilde f_I(u+d)\otimes \tilde f_J(u+e+p)$  we get precisely the right-hand side of (\ref{eq:yang-flag-YS-rels}), mutatis mutandis.

Similarly, using the very same propositions above, one can show that the tensor algebra $T(C_n)$ modulo the monomial straightening relations $(\mathcal M_{I,J})$ is a comodule algebra. It is left to check $(\mathcal C_{I,J})$, which we do in a more direct manner below.

\begin{mylemm} Let $\mathcal I$ be the ideal in $T(C_n)$ generated by the commuting relations $(\mathcal{C}_{I,J})$ of (\ref{eq:yang-flag-C-rels}). Then $\rho(\mathcal I) \subseteq I\otimes Y_n$, making the quotient $T(C_n) / \mathcal I$ a right $Y_n$-comodule algebra.
\end{mylemm}

\begin{proof} To save space, we replace the fraction $\frac{(-1)^{p-1} p!}{(u-v-e+1)\cdots (u-v-e+p)}$ depending on $p$ with $(\ast_p)$ throughout.

Applying $\rho$ to the left-hand side $(\star)$ of (\ref{eq:yang-flag-C-rels}), we have $\sum_{K,L} \tilde f_{K}(u) \tilde f_{L}(v) \otimes t^{K}_{J}(u) t^{L}_{I}(v) - \sum_{K,L} \tilde f_{L}(v) \tilde f_{K}(u)\otimes t^{L}_{I}(v) t^{K}_{J}(u)$.

Applying $\rho$ on the right $(\star\star)$ yields two terms as well,
$$\sum_{p=1}^{d} (\ast_p) \binom{d}{p}\sum_{K,L}\tilde f_{K}(u) \tilde f_{L}(v) \otimes t^{K}_{J}(u) t^{L}_{I}(v)
$$
and 
$$
-\sum_{p=1}^{d} (\ast_p) \sum_{A,B\in\binom{[n]}{p}} \sum_{K,L} \tilde f_{L}(v) \tilde f_{K}(u) \otimes t^{L}_{i_1\cdots j_{b_1}\cdots j_{b_p}\cdots i_d}(v) \cdot t^{K}_{j_1\cdots i_{a_1}\cdots i_{a_p}\cdots j_e}(u) .
$$
Rewrite the second term using (\ref{eq:yangian-C-rels}) and get
\begin{eqnarray*}
&&\sum_{K,L}\tilde f_{L}(v) \tilde f_{K}(u) \otimes \big[t^K_J(u)\,,\,t^L_I(v)\big] \\
&& - \sum_{K,L}\tilde f_{L}(v) \tilde f_{K}(u) \otimes \left(\sum_{p}(\ast_p)\sum_{A,B} t^{k_1\cdots \ell_{b_1}\cdots \ell_{b_p}\cdots k_e}_{J}(u) \cdot t^{\ell_1\cdots k_{a_1}\cdots k_{a_p}\cdots \ell_d}_{I}(v)\right).
\end{eqnarray*}
The first term of this last expression subtracts nicely from $\rho(\star)$:
\begin{eqnarray*}
&&\sum_{K,L} \left(\tilde f_K(u) \tilde f_L(v) - \tilde f_L(v) \tilde f_K(u) \right) \otimes t^K_J(u) t^L_I(v) \\
&& - \sum_{K,L} \left(\tilde f_L(v) \tilde f_K(u) - \tilde f_L(v) \tilde f_K(u) \right) \otimes t^L_I(v) t^K_J(u).
\end{eqnarray*}
Now use (\ref{eq:yang-flag-C-rels}) on the new left-hand side to get
\begin{eqnarray*}
&&\sum_{K,L} \sum_{p} (\ast_p)\binom{d}{p} \tilde f_K(u) \tilde f_L(v) \otimes t^K_J(u) t^L_I(v) \\
&& - \sum_{K,L}\left(\sum_{p}(\ast_p)\sum_{A,B} \tilde f_{\ell_1\cdots k_{a_1}\cdots k_{a_p}\cdots \ell_d}(v) \cdot \tilde f_{k_1\cdots \ell_{b_1}\cdots \ell_{b_p}\cdots k_e}(u) \right) \otimes t^K_J(u) t^L_I(v) 
\end{eqnarray*}
The first term here cancels the original first term of $\rho(\star\star)$, and we are left with demonstrating the equality of
\begin{equation}\label{eq:left-remains}
\sum_{p=1}^d (\ast_p)\sum_{K,L}\sum_{A,B} \tilde f_{\ell_1\cdots k_{a_1}\cdots k_{a_p}\cdots \ell_d}(v) \cdot \tilde f_{k_1\cdots \ell_{b_1}\cdots \ell_{b_p}\cdots k_e}(u) \otimes t^K_J(u) t^L_I(v) 
\end{equation}
and
\begin{equation}\label{eq:right-remains}
\sum_{p=1}^d (\ast_p)\sum_{K,L}\sum_{A,B} \tilde f_{L}(v) \tilde f_{K}(u) \otimes t^{k_1\cdots \ell_{b_1}\cdots \ell_{b_p}\cdots k_e}_{J}(u) \cdot t^{\ell_1\cdots k_{a_1}\cdots k_{a_p}\cdots \ell_d}_{I}(v)\,,
\end{equation}
which we may do one $p$-summand at a time. Note that, by the alternating property of $\tilde f_X$ and $t^Y_Z$, we may replace any summand in (\ref{eq:left-remains}) with
\begin{eqnarray*}
\nonumber && \sum_{K,L}\sum_{A,B}\left\{ (-1)^{\sum_{r} (a_r-r) }\tilde f_{K_A|L\setminus L_B}(v) \cdot (-1)^{\sum_{r}(b_r-r)} \tilde f_{L_B|K\setminus K_A}(u) \right\}\otimes \\
&&\left\{(-1)^{ \sum_{r} (a_r-r)}t^{K_A|K\setminus K_A}_J(u) \cdot (-1)^{\sum_{r}(b_r-r)} t^{L_B|L\setminus L_B}_I(v)\right\}\,.
\end{eqnarray*}
Similarly, a summand in (\ref{eq:right-remains}) reduces to 
\begin{equation*}
\sum_{K,L}\sum_{A,B} \tilde f_{L_B|L\setminus L_B}(v) \tilde f_{K_A|K\setminus K_A}(u) \otimes t^{L_B|K\setminus K_A}_J(u) t^{K_A|L\setminus L_B}_I(v)\,.
\end{equation*}
When $K_A\cap (L\setminus L_B)\neq \emptyset$ the summands involved above are zero. Likewise when $L_B\cap (K\setminus K_A)\neq \emptyset$. Let us denote this with Kronecker deltas. Also, we save space by denoting, e.g., $L\setminus L_B$ by $L^B$ and dropping the $u$'s and $v$'s. We must show the equality of 
\begin{equation*}
\sum_{K,L}\sum_{A,B}\left(\delta_{K_A,K^A}\delta_{K_A,L^B}\delta_{L_B,K^A}\delta_{L_B,L^B}\right)\tilde f_{K_A|L^B} \tilde f_{L_B|K^A} \otimes t^{K_A|K^A}_J t^{L_B|L^B}_I
\end{equation*}
and
\begin{equation*}
\sum_{K,L}\sum_{A,B} \left(\delta_{K_A,K^A}\delta_{K_A,L^B}\delta_{L_B,K^A}\delta_{L_B,L^B}\right)\tilde f_{L_B|L^B} \tilde f_{K_A|K^A} \otimes t^{L_B|K^A}_J t^{K_A|L^B}_I \,.
\end{equation*}
Now in our notation, $\delta_{K_A,K^A}$ and $\delta_{L_B,L^B}$ are obviously always $1$, but we include these because it allows us to rewrite the sum. Instead of summing over sets $K,L$ and then subsets $K_A, L_B$, let us some over sets $K_0, L_0$ and complements $K^+, L^+$. The previous two expressions become
\begin{equation*}
\sum_{K_0,L_0}\sum_{K^+,L^+}\left(\delta_{K_0,K^+}\delta_{K_0,L^+}\delta_{L_0,K^+}\delta_{L_0,L^+}\right)\tilde f_{K_0|L^+} \tilde f_{L_0|K^+} \otimes t^{K_0|K^+}_J t^{L_0|L^+}_I
\end{equation*}
and
\begin{equation*}
\sum_{K_0,L_0}\sum_{K^+,L^+}\left(\delta_{K_0,K^+}\delta_{K_0,L^+}\delta_{L_0,K^+}\delta_{L_0,L^+}\right)\tilde f_{L_0|L^+} \tilde f_{K_0|K^+} \otimes t^{L_0|K^+}_J t^{K_0|L^+}_I \,.
\end{equation*}
Finally, if we swap the labels $K_0$ and $L_0$ while leaving the labels $K^+$ and $L^+$ fixed
in the second expresion, we reach the first, concluding the proof of the lemma and the theorem.
\end{proof}

\subsection{Modules}
Here we return to the viewpoint that $Y_n$ is a deformation of $U(\mathfrak{gl}_n)$ and look for an action of $Y_n$ on $\mathcal{F}\ell_{Y_n}\!(\gamma)$ modeled after the classic setting.
For all $a,b\in[n]$ and any $J\in[n]^r$ with all entries distinct (though not necessarily arranged in order), define an action of $Y_n$ on $C_n(r)$ by
\begin{equation}\label{eq:yangian-action}
t_{ab}(u)\cdot \tilde f_{J}(v) = \delta_{ab} \tilde f_I(v) + \delta_{b\in J} u^{-1} \tilde f_{j_1\cdots a \cdots j_r}(v)\,.
\end{equation}
We show that this action: \emph{(i) is well-defined, i.e. it respects the relations (\ref{eq:def-rels-series}); (ii) extends to an action of $Y_n$ on $T(C_n)$; and (iii) preserves the ideal realizing $\mathcal{F}\ell_{Y_n}\!(\gamma)$ as a quotient of $T(C_n)$.} In other words, 

\begin{mytheo} $\mathcal{F}\ell_{Y_n}\!(\gamma)$ is a $Y_n$-module algebra.
\end{mytheo} 

\begin{proof}[Proof of i)] We must show that 
$$
[t_{ab}(u)\,,\,t_{cd}(v)]\cdot \tilde f_{J}(w) = \frac{1}{u-v}\Big(t_{cb}(u) t_{ad}(v) - t_{cb}(v) t_{ad}(u) \Big) \cdot \tilde f_J(w)\,,
$$
which we may break up into several cases: (1) $b=d$; (2) $b\neq d \wedge b=a$; (3) $b\neq d \wedge b\neq a \wedge b=c$; and (4) $b\neq d \wedge b\neq a \wedge b\neq c$. We skip the middle two cases for brevity and suppress the $w$'s for clarity. 
\smallskip

\noindent\textit{Case 1).}  On the left above we have
$$
\delta_{ab}\delta_{cb} \tilde f_J + \delta_{ab} \delta_{b\in J} \frac{1}{v} \tilde f_{j_1\cdots c\cdots j_r} + \delta_{cb} \delta_{b\in J} \frac{1}{u} \tilde f_{j_1\cdots a\cdots j_r} + \delta_{cb}\delta_{b\in J} \frac{1}{uv} \tilde f_{j_1\cdots a \cdots j_r}
$$
from $t_{ab}(u)t_{cb}(v)\cdot \tilde f_J$. The term $t_{cb}(v)t_{ab}(u)\cdot \tilde f_J$ looks similar, and after simplification, we have 
$$
u^{-1}v^{-1}\delta_{b\in J}\cdot\big( \delta_{cb} \tilde f_{j_1\cdots a \cdots j_r} - \delta_{ab} \tilde f_{j_1\cdots c\cdots j_r} \big)
$$
on the left-hand side.

On the right, we have $1/(u-v)$ times
$$
\delta_{cb} \delta_{b\in J} \frac{1}{v} \tilde f_{j_1\cdots a\cdots j_r} + \delta_{ab} \delta_{b\in J} \frac{1}{u} \tilde f_{j_1\cdots c\cdots j_r} - 
\delta_{cb} \delta_{b\in J} \frac{1}{u} \tilde f_{j_1\cdots a\cdots j_r} - \delta_{ab} \delta_{b\in J} \frac{1}{v} \tilde f_{j_1\cdots c\cdots j_r}  
$$
after similar simplifications; call this $(\star)$. Continuing, we have
\begin{eqnarray*}
(\star) &=& \frac{1}{u-v}\left(\delta_{cb} \delta_{b\in J} (v^{-1}-u^{-1}) \tilde f_{j_1\cdots a\cdots j_r} - \delta_{ab} \delta_{b\in J} (v^{-1}-u^{-1}) \tilde f_{j_1\cdots c\cdots j_r} \right)  \\
&=& u^{-1}v^{-1} \delta_{b\in J} \cdot \big(\delta_{cb} \tilde f_{j_1\cdots a \cdots j_r} - \delta_{ab} \tilde f_{j_1\cdots c\cdots j_r} \big),
\end{eqnarray*}
as needed.
\medskip

\noindent\textit{Case 4).} Here it will be useful to keep track of which symbol ($b$ or $d$) is being replaced. Let us augment our previous notation a bit: $t_{ab}(u) \cdot \tilde f_J =  u^{-1}\, \delta_{b\in J}\, \tilde f_{j_1 \cdots}\newstack{b}{a}{}_{\cdots j_r}$. Under the current hypotheses, the left-hand side becomes
\begin{align*}
\delta_{cd}\delta_{b\in J} \frac{1}{u} \tilde f_{j_1\cdots}\newstack{b}{a}{}_{\cdots j_r} &+ \delta_{b\in (J\setminus d)\cup c} \delta_{d\in J} \frac{1}{uv} \tilde f_{j_1\cdots}\newstack{b}{a}{}_{\cdots}\newstack{d}{c}{}_{\cdots j_r} \\
&- \delta_{cd}\delta_{b\in J} \frac{1}{u} \tilde f_{j_1\cdots}\newstack{b}{a}{}_{\cdots j_r} - \delta_{d\in (J\setminus b)\cup a} \delta_{b\in J} \frac{1}{uv} \tilde f_{j_1\cdots}\newstack{b}{a}{}_{\cdots}\newstack{d}{c}{}_{\cdots j_r}
\end{align*}
and the right-hand side becomes $\frac{1}{u-v}$ times
$$
t_{cb}(u) \left\{\delta_{ad} \tilde f_J + \delta_{d\in J} \frac{1}{v} \tilde f_{j_1\cdots}\newstack{d}{a}{}_{\cdots j_r} \right\} - t_{cb}(v) \left\{\delta_{ad} \tilde f_J + \delta_{d\in J} \frac{1}{u} \tilde f_{j_1\cdots}\newstack{d}{a}{}_{\cdots j_r} \right\}
$$
or $-u^{-1} v^{-1}\, \delta_{ad}\delta_{b\in J}\, \tilde f_{j_1\cdots}\newstack{b}{c}{}_{\cdots j_r}$.

Returning to the left-hand side, we notice that, under the hypotheses, $\delta_{b\in (J\setminus d)\cup c} = \delta_{b\in J}$, while $\delta_{d\in (J\setminus b)\cup a}$ acts as $\delta_{d\in J} + \delta_{da}$ on $\tilde f_{j_1\cdots}\newstack{b}{a}{}_{\cdots}\newstack{d}{c}{}_{\cdots j_r}$. We may omit the overlap case $a\in J$ because if this were true, the intermediate step $\tilde f_{j_1\cdots}\newstack{b}{a}{}_{\cdots j_r}$ would have produced a zero term. So we have
$$
\frac{1}{uv} \delta_{b\in J} \left(\delta_{d\in J} \tilde f_{j_1\cdots}\newstack{b}{a}{}_{\cdots}\newstack{d}{c}{}_{\cdots j_r} - \left\{ \delta_{d\in J} \tilde f_{j_1\cdots}\newstack{b}{a}{}_{\cdots}\newstack{d}{c}{}_{\cdots j_r} + \delta_{da} \tilde f_{j_1\cdots}\newstack{b}{a}{}_{\cdots}\newstack{d}{c}{}_{\cdots j_r} \right\}\right) ,
$$
or $-u^{-1} v^{-1}\, \delta_{ad}\delta_{b\in J}\, \tilde f_{j_1\cdots}\newstack{b}{\!\!\stackrel{a=d}{c}\!\!}{}_{\cdots j_r}=-u^{-1} v^{-1}\, \delta_{ad}\delta_{b\in J}\, \tilde f_{j_1\cdots}\newstack{b}{c}{}_{\cdots j_r}$, completing the proof in the final case.
\end{proof}

\begin{proof}[Proof of ii)] Given a monomial $\tilde f_{\vec{J}}(\vec{w}):= \tilde f_{J_1}(w_1) \tilde f_{J_2}(w_2) \cdots \tilde f_{J_p}(w_p)$ in $T(C_n)[[w_1^{-1}, \ldots, w_p^{-1}]]$, let us define an operator $\partial_{ab}^{\,i}$ for any $1\leq a,b,\leq n$ and $1\leq i\leq p$ as follows:
$$
\partial_{ab}^{\,i} \cdot \tilde f_{\vec{J}}(\vec w) = \tilde f_{J_1}(w_1) \cdots \left\{ \delta_{b\in J_i} \tilde f_{j_{i1}\cdots a \cdots j_{ir_i}}(w_i) \right\} \cdots \tilde f_{J_p}(w_p).
$$
Now define an action of $Y_n$ by
\begin{eqnarray*}
t_{ab}(u) \cdot \tilde f_{\vec{J}}(\vec{w}) &=& \delta_{ab} \tilde f_{\vec J}(\vec w) + u^{-1}\sum_{i=1}^{p}  \partial_{ab}^{\,i} \tilde f_{\vec{J}}(\vec w) 
\end{eqnarray*}
As before we drop the $w$'s appearing in the formulas to make the calculations more compact.
We must show that
$$
[t_{ab}(u)\,,\,t_{cd}(v)]\cdot \tilde f_{\vec J} = \frac{1}{u-v}\Big(t_{cb}(u) t_{ad}(v) - t_{cb}(v) t_{ad}(u) \Big) \cdot \tilde f_{\vec J}\,,
$$
the left-hand side of which is readily reduced to
\begin{equation*}
\sum_{k=1}^p \frac{1}{u}\partial_{ab}^{\,k} \sum_{i=1}^p \frac{1}{v}\partial_{cd}^{\,i} \tilde f_{\vec J} - \sum_{i=1}^p \frac{1}{v}\partial_{cd}^{\,i} \sum_{k=1}^p \frac{1}{u}\partial_{ab}^{\,k} \tilde f_{\vec J}\,.
\end{equation*}
The operators $\partial_{ab}^{\,k}$ and $\partial_{cd}^{\,i}$ commute when $i\neq k$. What remains on the left can be written as $\sum_{i} \tilde f_{J_1}\cdots \left\{\big[t_{ab}(u)\,,\, t_{cd}(v)\big] \cdot \tilde f_{J_i} \right\} \cdots \tilde f_{J_p}\,$, and a reduction to Part (i) looks likely.

From the right-hand side, we get $\frac{1}{u-v}$ times
\begin{equation*}
\delta_{ad} \sum_{i=1}^p \frac{1}{u}\partial_{cb}^{\,i} \tilde f_{\vec J} + \delta_{cb} \sum_{i=1}^p \frac{1}{v}\partial_{ad}^{\,i} \tilde f_{\vec J} - \delta_{ad} \sum_{i=1}^p \frac{1}{v}\partial_{cb}^{\,i} \tilde f_{\vec J} - \delta_{cb} \sum_{i=1}^p \frac{1}{u}\partial_{ad}^{\,i} \tilde f_{\vec J}\,,
\end{equation*}
or 
\begin{equation*}
\sum_{i=1}^p \tilde f_{J_1} \cdots \left\{\frac{1}{u-v} \left(t_{cb}(u)t_{ad}(v) - t_{cb}(v) t_{ad}(u)\right)\cdot \tilde f_{J_i} \right\}\cdots \tilde f_{J_p}\,;
\end{equation*}
confirming our suspicions about Part (i).
\end{proof}

\begin{proof}[Proof of iii)] One must check that the action respects the alternating, Young symmetry, monomial straightening, and commuting relations. The first check is easy and the third looks much like the second, so we omit them.
\bigskip

\noindent\textit{Proof of }$(\mathcal Y_{I,J})$: \newline
Fix $d\leq e\in\Vert\gamma\Vert$, $1\leq r$, $I\in\binom{[n]}{d}$, and $J\in\binom{[n]}{e}$. We show that
$$
t_{ab}(u) \cdot \sum_{\Lambda\in\binom{[n]}{r}} (-1)^{\Ell{J_{\Lambda}}{J\setminus J_{\Lambda}}} \tilde f_{I | J_{\Lambda}}(v+\alpha) \tilde f_{J\setminus J_{\Lambda}}(v+\beta) \equiv 0
$$
modulo the ideal in $T(C_n)$ generated by the young symmetry relations. Writing out the definition of the action, straightaway we are left with showing that
\begin{eqnarray*}
&&\sum_{\Lambda\in\binom{[n]}{r}} (-1)^{\Ell{J_{\Lambda}}{J\setminus J_{\Lambda}}} \partial_{ab} \tilde f_{I | J_{\Lambda}}(v+\alpha) \tilde f_{J\setminus J_{\Lambda}}(v+\beta) \\
&&+ \sum_{\Lambda\in\binom{[n]}{r}} (-1)^{\Ell{J_{\Lambda}}{J\setminus J_{\Lambda}}} \tilde f_{I | J_{\Lambda}}(v+\alpha) \partial_{ab} \tilde f_{J\setminus J_{\Lambda}}(v+\beta) 
\end{eqnarray*}
is congruent to zero. Now, the first involves the Kronecker delta function $\delta_{b\in I|J_{\Lambda}}$, which we first write as $\delta_{b\in I} + \delta_{b\in J_{\Lambda}} - \delta_{b\in I\cap J_{\Lambda}}$. Of course, if $I\cap J_{\Lambda}$ is ever nonempty, then $t_{ab}(u)$ will never see the corresponding summand because $\tilde f_{I|J_{\Lambda}}=0$. The function $\delta_{b\in I}$ shows up above as
$$
\delta_{b\in I} \sum_{\Lambda} (-1)^{\Ell{J_{\Lambda}}{J\setminus J_{\Lambda}}} \tilde f_{i_1\cdots a \cdots i_{d-r} | J_{\Lambda}}(v+\alpha) \tilde f_{J\setminus J_{\Lambda}}(v+\beta) \,,
$$
which is another Young symmetry relation, hence congruent to zero. We are left with
\begin{eqnarray*}
&&\sum_{\Lambda\in\binom{[n]}{r}} (-1)^{\Ell{J_{\Lambda}}{J\setminus J_{\Lambda}}} \times \\
&& \left\{\delta_{b\in J_{\Lambda}} \tilde f_{I | j_{\lambda_1}\cdots a \cdots j_{\lambda_r}}(v+\alpha) \tilde f_{J\setminus J_{\Lambda}}(v+\beta)  + \tilde f_{I | J_{\Lambda}}(v+\alpha) \partial_{ab} \tilde f_{J\setminus J_{\Lambda}}(v+\beta) \right\},
\end{eqnarray*}
only one term of which is nonzero for any given $\Lambda$. We may rewrite this sum as
$$
\delta_{b\in J} \sum_{\Lambda\in\binom{[n]}{r}} (-1)^{\Ell{J'_{\Lambda}}{J'\setminus J'_{\Lambda}}} \tilde f_{I | J'_{\Lambda}}(v+\alpha) \tilde f_{J'\setminus J'_{\Lambda}}(v+\beta)\,,
$$
where if $J = (j_1, \ldots, b, \ldots, j_{e+r})$, then $J' = (j_1,\ldots, a, \ldots j_{e+r})$; this is another Young symmetry relation.
\medskip

\noindent\textit{Proof of }$(C_{I,J})$: \newline
Fix $d\leq e\in\Vert\gamma\Vert$, $I\in[n]^d$, and $J\in[n]^e$. The expression $t_{ab}(u)\cdot \big[\tilde f_{J}(v)\,,\, \tilde f_{I}(w) \big]$ simplifies to
$$
\delta_{ab} \big[\tilde f_{J}(v)\,,\, \tilde f_{I}(w)\big] + \delta_{b\in J} \frac{1}{u}\big[\tilde f_{j_1\cdots a \cdots j_e}(v)\,,\, \tilde f_{I}(w)\big] + \delta_{b\in I} \frac{1}{u}\big[\tilde f_{J}(v)\,,\, \tilde f_{i_1\cdots a\cdots i_d}(w)\big].
$$
The above should be the same as $t_{ab}(u)$ applied to
$$
\sum_{p=1}^d (\ast_p)  \bigg\{ \binom{d}{p} \tilde f_{J}(v) \tilde f_{I}(w) - \sum_{K,L\in\binom{[n]}{p}} \tilde f_{i_1\cdots j_{\ell_1}\cdots j_{\ell_p}\cdots i_d}(w) \tilde f_{j_1\cdots i_{k_1} \cdots i_{k_p} \cdots j_e}(v) \bigg\},
$$
Let us simplify notation a bit. First, drop the $v$'s and $w$'s appearing here. Second, write, e.g., $(i_1, \ldots, j_{\ell_1}, \ldots, j_{\ell_p}, \ldots, i_d)$ as $I^K\!{}_{\curlywedge} J_L$. We get
\begin{eqnarray*}
&& \delta_{ab}\sum_{p=1}^d (\ast_p)  \bigg\{ \binom{d}{p} \tilde f_{J} \tilde f_{I} - \sum_{K,L\in\binom{[n]}{p}} \tilde f_{I^K\!{}_{\curlywedge} J_L} \tilde f_{J^L\!\!{}_{\curlywedge} I_K} \bigg\} \\
&& +\frac{1}{u}\sum_{p} (\ast_p) \binom{d}{p} \bigg\{ \partial_{ab} \tilde f_J \tilde f_I + \tilde f_J  \partial_{ab} \tilde f_I \bigg\} \\
&& -\frac{1}{u}\sum_{p} (\ast_p) \sum_{K,L} \bigg\{ \partial_{ab} \tilde f_{I^K\!{}_{\curlywedge} J_L} \tilde f_{J^L\!\!{}_{\curlywedge} I_K} + \tilde f_{I^K\!{}_{\curlywedge} J_L} \partial_{ab} \tilde f_{J^L\!\!{}_{\curlywedge} I_K} \bigg\} .
\end{eqnarray*}
Notice that the function $\delta_{b\in (I^K\!{}_{\curlywedge} J_L)}$ appearing in the term $\partial_{ab}\tilde f_{I^K\!{}_{\curlywedge} J_L} \tilde f_{J^L\!\!{}_{\curlywedge} I_K}$ above takes the same value as $\delta_{b\in I\setminus I_K} + \delta_{b\in J_L}$, since any summand satisfying $(I\setminus I_K)\cap J_L\neq\emptyset$ vanishes. Similarly rewriting the function $\delta_{b\in (J^L\!\!{}_{\curlywedge} I_K)}$ appearing in $\tilde f_{I^K\!{}_{\curlywedge} J_L} \partial_{ab} \tilde f_{J^L\!\!{}_{\curlywedge} I_K}$ and rearranging the sums, we may write the above as
\begin{eqnarray*}
&& \delta_{ab}\sum_{p=1}^d (\ast_p)  \bigg\{ \binom{d}{p} \tilde f_{J} \tilde f_{I} - \sum_{K,L\in\binom{[n]}{p}} \tilde f_{I^K\!{}_{\curlywedge} J_L} \tilde f_{J^L\!\!{}_{\curlywedge} I_K} \bigg\} \\
&& +\delta_{b\in J} \frac{1}{u}\sum_{p} (\ast_p) \bigg\{ \binom{d}{p} \tilde f_{J'} \tilde f_I   - \sum_{K,L} \tilde f_{I^K\!{}_{\curlywedge} J'_L} \tilde f_{{J'}^L\!\!{}_{\curlywedge} I_K} \bigg\} \\
&& + \delta_{b\in I} \frac{1}{u}\sum_{p} (\ast_p) \bigg\{ \binom{d}{p} \tilde f_J \tilde f_{I'} -\sum_{K,L} \tilde f_{{I'}^K\!{}_{\curlywedge} J_L} \tilde f_{J^L\!\!{}_{\curlywedge} I'_K} \bigg\},
\end{eqnarray*}
where again, e.g., $J' = (j_1, \ldots, a, \ldots j_{e})$ when $J=(j_1, \ldots, b, \ldots, j_e)$. Compare this to what we had on the left---in the new notation
$$
\delta_{ab} \big[\tilde f_{J}\,,\, \tilde f_{I}\big] + \delta_{b\in J} \frac{1}{u}\big[\tilde f_{J'}\,,\, \tilde f_{I}\big] + \delta_{b\in I} \frac{1}{u}\big[\tilde f_{J}\,,\, \tilde f_{I'}\big].
$$
Conclude the two sides agree modulo the ideal generated by $(C_{I,J})$.
\end{proof}

Note that $t_{ii}(u) f_{\vec{J}}(v)$ is a $\C[u^{-1}]$-multiple of $f_{\vec{J}}(v)$ for all set-tuples $\vec{J}$ and all $1\leq i\leq n$. Moreover, $t_{ij}(u)\cdot f_{\vec{J}}(v)=0$ whenever $i<j$ and $\vec{J}=([d_1],\ldots, [d_p])$. The reader has by now guessed that a Yangian version of the highest-weight theory in Section \ref{sec:comm-modules} is known to hold, cf. \cite{Dri:1,BilFutMol:1}. In our setup, it is not immediately clear what generator $v$, if any, satisfies $Y_n\cdot v = \mathcal{F}\ell_{Y_n}\!(\gamma)$. However, the analog of Theorem \ref{thm:comm-basis-thm} holds---the preferred basis monomials are $f_{I_1}^{r_1}\cdots f_{I_p}^{r_p}$ with the sizes of the $I_k$ now increasing, and with the $r_k$ arbitrary. We expect an analog of Theorem \ref{thm:comm-module-model} to hold as well.

\subsection{Parabolic Presentations}
Returning to the quasi-Pl\"ucker coordinates, we make a connection between noncommutative flags and the parabolic presentations of $Y_n$ given by Brundan and Kleshchev \cite{BruKle:1}. 
In the proof of Theorem \ref{thm:g-invt-fcns}, we factor the matrix $A=(a_{ij})$ as $\mathbb L\cdot \mathbb D\cdot \mathbb U$ inside the field $F\lskew A \rskew$.\footnote{We did not make $\mathbb L$ or $\mathbb D$ explicit, but they are filled with right/row quasi-Pl\"ucker coordinates and the quasideterminants $|A_{[d],[d]}|_{dd}$ respectively.} Notice that up to Equation (\ref{eq:LDU-std}), the only divisions carried out in the factorization are by elements $|A_{[d],[d]}|_{dd}$. Letting $A$ be the matrix of generators $T(u)$ for $Y_n$, this means (\ref{eq:LDU-std}) may be reached entirely within $Y_n[[u^{-1}]]$, with no need to pass to the larger skew field $D$ to carry out the calculations. For on the one hand, the series $t^{[d]}_{[d]}(u)$ starts with $1$ and may be inverted in $Y_n[[u^{-1}]]$, while on the other hand,  $|T_{[d],[d]}(u)|_{dd}$ is just $t^{[d]}_{[d]}(u+d-1)\cdot {t^{[d-1]}_{[d-1]}(u+d-1)}^{-1}$ by (\ref{eq:ydet-qdet-factors}).  

Brundan and Kleshchev show that: (i) the nonzero entries of $\mathbb L$, $\mathbb D$, and $\mathbb U$ all belong to $Y_n[[u^{-1}]]$, not just $D$ (just reverified above); and (ii) the subalgebra generated by $\mathcal G$ the set of coefficients of the powers of $u^{-1}$ appearing in the nonzero entries of $\mathbb L$, $\mathbb D$, and $\mathbb U$ actually generate all of $Y_n$ (obvious after \emph{unfactoring}, e.g. noting that $t_{ij}(u)$ is just the sum $\sum_{k} \mathbb L_{ik} \mathbb D_{kk} \mathbb U_{kj}$). The nontrivial part of \cite{BruKle:1} is as follows: they describe relations $\mathcal R$ among the generators $\mathcal G$ and show that these are a necessary and sufficient to define $Y_n$ abstractly as $\C\langle\mathcal G\rangle / \mathcal R$.

After (\ref{eq:LDU-std}), we have another description of these generators. Fix a composition $\gamma \models n$ as usual, and let $d_1, d_2, \ldots, d_r$ again denote the partial sums $\Vert\gamma\Vert$. For any $1\leq a <r$, choose $i,j$ satisfying $d_a < i\leq d_{a+1}$ and $d_{a+1} < j < n$. Then the $(i,j)$-entry of $\mathbb L$ is simply $p_{ij}^{[d_a]\setminus i}(T(u))$. In other words, the generators in $\mathcal G$ coming from $\mathbb L$ are just the coefficients of the powers of $u^{-1}$ occuring in the (left) quasi-Pl\"ucker coordinates.\footnote{Similarly, the generators coming from $\mathbb U$ are related to right/row quasi-Pl\"ucker coordinates.} It is interesting to note that the relation 
$$
1 = \sum_{1\leq j\leq e} p_{ij}^{[d]\setminus i} \cdot p_{ji}^{[e]\setminus j}.
$$
does not appear in $\mathcal{R}$. Also noteworthy: among all the quasi-Plucker coordinate relations given in Section \ref{sec:qplucker-coords}, this is the only one holding inside $Y_n[[u^{-1}]]$ (i.e., not requring the skew field $D$ to describe). It may be interesting to see how and if relations of this type can simplify $\mathcal R$. In the reverse direction, one might find new quasi-Pl\"ucker relations for $T(u)$ (and ideally, for generic matrices $A$) by studying the set $\mathcal R$. 

\def\cprime{$'$}
\providecommand{\bysame}{\leavevmode\hbox to3em{\hrulefill}\thinspace}
\providecommand{\MR}{\relax\ifhmode\unskip\space\fi MR }
\providecommand{\MRhref}[2]{%
  \href{http://www.ams.org/mathscinet-getitem?mr=#1}{#2}
}
\providecommand{\href}[2]{#2}

\bigskip

\par\noindent
{\small\textsc{LaCIM, UQ\`aM, Case Postale 8888,
succursale Centre-ville, Montr\'eal (Qu\'ebec) H3C 3P8, Canada.}\\
\texttt{lauve@lacim.uqam.ca}}

\end{document}